\documentclass [10pt,reqno]{amsart}
\usepackage{ucs}
\usepackage{amsmath, amssymb, amscd}
\usepackage{pb-diagram}
\usepackage[latin1]{inputenc}
\usepackage{amsmath}
\usepackage{amssymb}
\usepackage{epsfig}
\usepackage{graphicx}
\usepackage{psfrag}

\newtheorem{theorem}{Theorem}[section]
\newtheorem{definition}[theorem]{Definition}
\newtheorem{lemma}[theorem]{Lemma}

\newtheorem{prop}[theorem]{Proposition}
\newtheorem{corollary}[theorem]{Corollary}

\newtheorem*{observation}{Observation}
\newtheorem*{remark}{Remark}

\newcommand{\norm}[1]{\| #1 \|}

\renewcommand{\epsilon}{\varepsilon}

\DeclareMathAlphabet{\mathpzc}{OT1}{pzc}{m}{it}
\usepackage{mathrsfs}

\newcommand{\Z}{\mathbb{Z}}
\newcommand{\C}{\mathbb{C}}

\newcommand{\R}{\mathbb{R}}

\renewcommand{\qed}{$\hfill \square$ \smallskip \\}
\renewcommand{\phi}{\varphi}

\newcommand{\ad}{\text{ad}}
\newcommand{\Ad}{\text{Ad}}

\newcommand{\Aut}{\text{Aut}}
\renewcommand{\det}{\text{det}}

\newcommand{\su}{\mathfrak{su}}

\newcommand{\G}{\mathscr{G}}

\begin{document}
\thispagestyle{empty}
\title[Casson-type moduli spaces over definite 4-manifolds]{On Casson-type instanton moduli spaces over
negative definite four-manifolds}
\author{Andrew Lobb \\ Raphael Zentner}

\begin {abstract} 
Recently Andrei Teleman considered instanton moduli spaces over negative definite four-manifolds $X$ with
$b_2(X) \geq 1$. If $b_2(X)$ is divisible by four and $b_1(X) =1$ a gauge-theoretic invariant can be defined;
it is a count of flat connections modulo the gauge group. Our first result shows that if such a moduli space
is non-empty and the manifold admits a connected sum decomposition $X \cong X_1 \# X_2$ then both $b_2(X_1)$
and $b_2(X_2)$ are divisible by four; this rules out a previously naturally appearing source of $4$-manifolds
with non-empty moduli space.  We give in some detail a construction of negative definite $4$-manifolds which
we expect will eventually provide examples of manifolds with non-empty moduli space.
\end {abstract}

\address{Mathematics Department \\ Imperial College London \\ London SW11 7AZ \\ UK}
\email{a.lobb@imperial.ac.uk}
\address {Fakult\"at f\"ur Mathematik \\ Universit\"at Bielefeld \\ 33501 Bielefeld\\
Germany}
\email{rzentner@math.uni-bielefeld.de}

\maketitle

\section*{Introduction}
Recently Andrei Teleman considered moduli spaces of projectively anti-selfdual instantons in
certain Hermitian rank-2 bundles over a closed oriented 4-manifold with negative definite
intersection form \cite{T}. These play a role in his classification program on Class VII surfaces \cite{T2}\cite{T3}. However, in
certain situations the instanton moduli spaces involved consist of projectively flat connections and
therefore have very interesting topological implications. In this article we will study these `Casson-type'
moduli spaces. 

 Suppose $E \to X$ is a Hermitian rank-2 bundle with first Chern-class a (minimal)
characteristic vector $w$ of the intersection form. In other words, it is the sum of elements $\{ e_i \}$ in
$H^2(X;\Z)$ which induce a basis of
$H^2(X;\Z)/\text{\em Tors}$ diagonalising the intersection form (because of Donaldson's theorem \cite{D}). 
Then for one possible value of a strictly negative second Chern class $c_2(E)$ the moduli space is compact
(independently of the Riemannian metric). In particular, if the manifold has second Betti-number $b_2(X)$
divisible by 4 and first Betti-number $b_1(X) = 1$ the instanton moduli space consists of
projectively flat connections and has expected dimension zero.  This should be thought of as a `Casson-type' moduli space because the holonomy yields a surjection onto the space of
$SO(3)$ representations of $\pi_1(X)$ with fixed Stiefel-Whitney class $w_2 = w \ (mod \ 2)$.

Non-emptiness of the Casson-type moduli space implies that none of the
elements $e_i$ can be Poincar\'e dual to an element
representable by a sphere, i.e. to an element in the image of the Hurewicz homomorphism. Prasad and Yeung
\cite{PY} constructed aspherical manifolds $W$ which are rational-cohomology complex projective planes, generalisations of Mumford's fake projective plane \cite{M}. If $\overline{W}$ denotes this manifold with the opposite orientation, a natural candidate of a manifold for which the moduli space might be non-empty is given
by the connected sum $4\overline{W}$ of 4 copies of $\overline{W}$, and a candidate of a manifold for which
the Casson-invariant can be defined is given by a `ring of 4 copies of $\overline{W}$' (the last summand in the 4-fold connected sum $4 \overline{W}$ is taken a connected sum with the first).

After recalling the gauge-theoretical situation considered in \cite{T} we show that if the Casson-type
moduli space is non-empty, then we cannot have a connected sum
decomposition $X \cong X_1 \# X_2$ unless both $b_2(X_1)$ and $b_2(X_2)$ are divisible by four. In
particular the moduli space for the above mentioned $4\overline{W}$ - ring is empty. 

This result still leaves open the question of whether there is any $X$ with a non-empty Casson-type
moduli space.  We give therefore in some detail a possible construction of suitable $4$-manifolds $X$ (along
with the correct representations of $\pi_1(X)$). We would like to point out that even though recent investigation leads us to believe that the Casson-type invariant is vanishing \cite{Z}, the Casson-type moduli space may still be non-empty and is interesting from a topological perspective. Our construction also suggests the possibility of considering Casson-type moduli spaces for manifolds with boundary.

\begin{remark}
A similar moduli space and invariant has been defined by Ruberman and Saveliev for $\Z[\Z]$-homology Hopf surfaces, going
back to work of Furuta and Ohta \cite{FO}, and for $\Z[\Z]$-homology 4-tori \cite{RS}. Our situation is
simpler than their first mentioned situation because of the absence of reducibles in the moduli space due to
the condition on $b_2(X)$.
\end{remark}

\section*{Acknowledgements}
The first author thanks Simon Donaldson for useful conversations.  The second author is grateful to Andrei
Teleman for turning his interest to low-energy instantons and for a stimulating conversation on them, and
also wishes to express his gratitude to Stefan Bauer for helpful conversations. Both authors thank Kim
Fr\o yshov profusely for invaluable advice and ideas.  We are also grateful to the referee for the care taken in helping us substantially improve the article.

\section{Donaldson theory on negative definite four-manifolds, low-energy instantons}
After briefly recalling some general instanton gauge theory \cite{DK}, and introducing our notations, we
shall quickly turn to the special situation of `low-energy instantons' over negative definite 4-manifolds
mentioned in the introduction. We show that the gauge-theoretical situation is indeed relatively simple,
indicate a definition of an invariant, and set up the correspondance of the moduli space to representation
spaces of the fundamental group in $SO(3)$.

\subsection{Connections}
Let $X$ be a smooth Riemannian four-manifold and $E \to X$ a Hermitian rank-2 bundle on $X$. Let further $a$
be a fixed unitary connection in the associated determinant line bundle $det(E)\to X$. We define
$\mathscr{A}_a(E)$ to be the affine space of unitary connections on $E$ which induce the fixed connection
$a$ in $det(E)$. This is an affine space over $\Omega^1(X;\su(E))$, the vector space of $\su(E)$-valued
one-forms on $X$. Let us denote by $P$ the principal $U(2)$ bundle of frames in $E$, and let $\overline{P}$
be the bundle that is associated to $P$ via the projection $\pi: U(2) \to PU(2)$, $\overline{P}=P \times_\pi
PU(2)$. The space $\mathscr{A}(\overline{P})$ of connections in the $PU(2)$ principal bundle $\overline{P}$
and the space $\mathscr{A}_a(E)$ are naturally isomorphic. If we interpret a connection $A \in
\mathscr{A}_a(E)$ as a $PU(2)$ connection via this isomorphism it is commonly called a projective
connection. The adjoint representation $\ad: SU(2) \to SO(\su(2))$ descends to a group isomorphim $PU(2) \to
SO(\su(2))$. The associated real rank-3 bundle $\overline{P} \times_{\ad} \su(2)$ is just the bundle
$\su(E)$ of traceless skew-symmetric endomorphisms of $E$. Thus the space $\mathscr{A}_a(E)$ is also
isomorphic to the space $\mathscr{A}(\su(E))$ of linear connections in $\su(E)$ compatible with the metric.
We shall write $A \in \mathscr{A}(\overline{P})$ for connections in the $PU(2)$ principal bundle and denote
the associated connection in $\su(E)$ by the same symbol. Should we mean the unitary connection which
induces the connection $a$ in $\det(E)$ we will write $A_a$ instead. 

Let $\G^0$ denote the group of automorphisms of $E$ of determinant $1$. It is called the `gauge group'. This
group equals the group of sections $\Gamma(X;P \times_\Ad SU(2))$, where $\Ad: U(2) \to \Aut(SU(2))$ is
given by conjugation. We shall write $\mathscr{B}(E)$ for the quotient space $\mathscr{A}(\overline{P}_E) /
\G^0$. A connection is called {\em reducible} if its stabiliser under the gauge group action equals the
subgroup given by the centre $\Z/2 = Z(SU(2))$ which always operates trivially, otherwise {\em irreducible}.
Equivalently, a connection $A_a$ is reducible if and only if there is a $A_a$ - parallel splitting of $E$
into two proper subbundles.

Let us point out that the characteristic classes of the bundle $\su(E)$ are given by
\begin{equation}\label{char classes}
\begin{split}
 	w_2(\su(E)) & = c_1(E) \ (mod \ 2) \\
	p_1(\su(E)) & = -4 c_2(E) + c_1(E)^2 \ .
\end{split}
\end{equation}

\subsection{Moduli space of anti-selfdual connections}
For a connection $A \in \mathscr{A}(\overline{P})$ we consider the anti-selfduality equation
\begin{equation}\label{asd}
F_A^+ = 0 \ ,
\end{equation}
where $F_A$ denotes the curvature form of the connection $A$, and $F_A^+$ its self-dual part with respect to
the Hodge-star operator defined by the Riemannian metric on $X$. The moduli space $\mathscr{M}(E) \subseteq
\mathscr{B}(E)$ of antiself-dual connections,
\begin{equation*}
\mathscr{M}(E) = \left.  \{ A \in \mathscr{A}(\overline{P}_E) \right| F_A^+ = 0 \} / \G^0
\end{equation*}
is the central object of study in instanton gauge theory. This space is in general non-compact and there is
a canonical ``Uhlenbeck-compactification'' of it. The anti-selfduality equations are elliptic, so Fredholm
theory provides finite dimensional local models for the moduli space. The often problematic aspect of
Donaldson theory is the need to deal with reducible connections and with a non-trivial compactification. We
will consider special situations where these problems do not occur.

\subsection{Low-energy instantons over negative definite four-manifolds}
We restrict now our attention to smooth Riemannian four-manifolds $X$ with $b_2^+(X) = 0$ and $b_2(X) \geq
1$.  According to Donaldson's theorem \cite{D} the intersection form of such a four-manifold is diagonal.  Let $\{e_i\}$ be a set of elements in $H^2(X;\Z)$ which induce a
basis of $H^2(X;\Z)/\text{\em Tors}$ diagonalising the intersection form.
\begin{lemma}\label{no reductions} \cite[section 4.2.1]{T}
Suppose the Hermitian rank-2 bundle $E \to X$ has first Chern class $c_1(E) = \sum e_i$ and its second Chern
class is strictly negative, $c_2(E) < 0$. Then $E \to X$ does not admit any topological decomposition $E = L
\oplus K$ into the sum of two complex line bundles. 
\end{lemma}
{\em Proof:} Suppose $E = L \oplus K$. Then $c_1(L) = \sum l_i \, e_i$ and $c_1(K)= \sum e_i - \sum l_i e_i$ for some $l_i \in \Z$. Therefore, 
\begin{equation*}
c_2(E) = c_1(L)(c_1(E) - c_1(L)) = \sum (l_i^2 - l_i) \geq 0 \ .
\end{equation*}
\hfill $\square$
\begin{corollary}
Let $E \to X$ be as in the previous lemma. Then the moduli space $\mathscr{M}(E)$ does not admit
reducibles.
\end{corollary}

For a connection $A \in \mathscr{A}(\su(E))$ Chern-Weil theory gives the following formula:
\begin{equation}\label{chern-weil formula}
\begin{split}
\frac{1}{8\pi^2} (\norm{F_A^-}^2_{L^2(X)} - \norm{F_A^+}^2_{L^2(X)})  = \, - \frac{1}{4} \, p_1(\su(E)) 
	 = c_2(E) - \frac{1}{4} c_1(E)^2
\end{split}
\end{equation}
In particular, for anti-selfdual connections the left hand side of this equation is always non-negative, and we
can draw the following observation from the formula:
\begin{observation} \cite[p. 1717]{T}
1. For $c_2(E) - 1/4 \, c_1(E)^2 \in \{0,1/4,2/4,3/4\}$ the moduli space $\mathscr{M}(E)$ is always compact,
independently of the chosen metric or any genericity argument. In fact, the lower strata in the
Uhlenbeck-compactification consist of anti-selfdual connections in bundles $E_k$ with $c_1(E_k) = c_1(E)$
and $c_2(E_k) = c_2(E) - k$ for $k \geq 1$. 
\\
2. For $c_1(E) = \sum e_i$ we have $c_1(E)^2 = -b_2(X)$. Thus, if $b_2(X) \equiv 0 \ (mod \ 4)$ and $c_2(E)
= - \frac{1}{4} \, b_2(X)$ the moduli space $\mathscr{M}(E)$ will consist of projectively flat connections
only.
\end{observation}

We recall the expected dimension of the moduli space $\mathscr{M}(E)$. It is given by the formula
\begin{equation*}
d(E) = -2 \, p_1(\su(E)) + 3  (b_1(X) - b_2^+(X) -1) 
\end{equation*}
In particular it can happen that $d(E) \geq 0$ in the situation we consider, namely, $b_2^+(X) = 0$, $c_1(E)
= \sum e_i$, and $c_2(E) < 0$, the latter condition assuring that we are in the favorable situation of Lemma
\ref{no reductions}.
\\

Interesting is the following special case of `Casson-type' moduli spaces that we consider from now on:
\begin{prop}\label{flat}
Let $X$ be a negative definite Riemannian four-manifold with strictly positive second Betti-number $b_2(X)$
divisible by four, and $b_1(X) = 1$. Let $E \to X$ be a Hermitian rank-2 bundle with $c_1(E) = \sum e_i$
and with $c_2(E) = -1/4 \ b_2(X)$. Then the moduli space $\mathscr{M}(E)$ of projectively anti-selfdual
connections in $E$ is compact and consists of irreducible projectively flat connections only, and is of
expected dimension zero.
\end{prop}

After suitable perturbations a gauge theoretic invariant can be defined in this situation: It is an algebraic count of a perturbed moduli space which consists of a finite number of points, the sign of each point is obtained by a natural orientation determined by the determinant line bundle of a family of elliptic operators. This has been done in the meantime in \cite{Z}, where it is shown that this invariant is actually zero. We would like to emphasise that the vanishing of this invariant doesn't imply emptiness of the unperturbed moduli space that we shall investigate further here.

\subsection{Flat connections, holonomy and representations of the fundamental group}
Suppose we are in the situation that our moduli space $\mathscr{M}(E)$ consists of flat connections in
$\su(E) \to X$, as for instance in the last proposition. Then we must have $p_1(\su(E)) = 0$ by Chern-Weil
theory. 

The holonomy establishes a correspondance between flat connections in the oriented real rank-3 bundle $V \to
X$ and representations of the fundamental group $\pi_1(X)$ in $SO(3)$ with a prescribed Stiefel-Whitney
class. More precisely, let $\rho : \pi_1(X) \to SO(3)$ be a representation of the fundamental group. Let
$\widetilde{X}$ be the universal covering of $X$; it is a $\pi_1(X)$ principal bundle over $X$. We can form
the associated oriented rank-3-bundle
\[
V_\rho := \widetilde{X} \times_\rho \R^3 \ .
\]
It admits a flat connection as it is a bundle associated to a principal bundle with discrete structure
group. Therefore it has vanishing first Pontryagin class, $p_1(V_\rho) = 0$, by Chern-Weil theory. Its only
other characteristic class \cite{DW} is its second Stiefel-Whitney class $w_2(V_\rho)$. Therefore we will
say that the representation $\rho$ has Stiefel-Whitney class $w \in H^2(X;\Z/2)$ if $w = w_2(V_\rho)$. On
the other hand, let $V \to X$ be an oriented real rank-3 bundle with a flat connection $A$. Then the
holonomy of $A$ along a path only depends up to homotopy on the path, and therefore induces a representation
$\text{\em Hol}(A) : \pi_1(X) \to SO(3)=SO(V|_{x_0})$. In particular, the holonomy defines a reduction of
the structure group to $\pi_1(X)$, and the bundle can therefore be reconstructed as $V \cong V_{\text{\em
Hol}(A)}$. In particular the representation $\text{\em Hol}(A)$ has Stiefel-Whitney class $w_2(V_{\text{\em
Hol}(A)}) = w_2(V)$. 

The moduli space $\mathscr{M}(E)$ has been obtained by quotienting the space of antiself-dual connections in
$\mathscr{A}(P_E) \cong \mathscr{A}(\su(E))$ by the gauge group $\G^0$. From the perspective of the $PU(2)$
connections in $\su(E)$ this gauge group is not the most natural one. Instead, the group
\[
\mathscr{G} := \Gamma(X;P \times_\Ad PU(2))
\]
is the natural group of automorphisms of connections in $\su(E)$. Not every element $g \in \G$ admits a lift
to $\G^0$; instead, there is a natural exact sequence
\begin{equation*}
1 \to \G^0 \to \G \to H^1(X;\Z/2) \to 0 \ .
\end{equation*}
Quotienting by $\mathscr{G}^0$ has the advantage of a simpler discussion of reducibles, as discussed above.
Let us denote by 
\begin{equation*}
\mathscr{M}(\su(E)):= \{A \in \mathscr{A}(\su(E)) \, | \, F_A^+ = 0 \} / \G \ 
\end{equation*}
the moduli space of anti-self dual connections in $\su(E)$ modulo the full gauge group $\G$. Then there is a
branched covering $\mathscr{M}(E) \to \mathscr{M}(\su(E))$ with `covering group' $H^1(X;\Z/2)$. 

Let us denote by $\mathscr{R}_w(\pi_1(X);SO(3))$ the space of representations of $\pi_1(X)$ in $SO(3)$ up to
conjugation and of Stiefel-Whitney class $w \in H^2(X;\Z/2)$. 
The above discussion implies that there is a homeomorphism
\[
\text{\em Hol}: \mathscr{M}(\su(E)) \stackrel{\cong}{\to} \mathscr{R}_w(\pi_1(X);SO(3)) \ ,
\]
where $w = w_2(\su(E))$. In particular, $\mathscr{M}(E)$ surjects onto $\mathscr{R}_w(\pi_1(X);SO(3))$.

\section{Representations of the fundamental group in $SO(3)$ and the vanishing result}

We will use the above derived relation of the `Casson-type moduli space' $\mathscr{M}(E)$ to the
representation space $\mathscr{R}_w(\pi_1(X);SO(3))$ to obtain the vanishing result which is mentioned in
the introduction. 

\subsection{Flat $SO(3)$ bundles}

The above construction of the bundle $V_\rho$ associated to a representation $\rho: \pi_1(X) \to SO(3)$ is
functorial in the following sense:
\begin{lemma}\label{naturality}
Suppose we have a map $f: W \to X$ between topological spaces, and $\rho : \pi_1(X) \to SO(3)$ a
representation of the fundamental group of $X$. Then there is a natural isomorphism
\begin{equation}
f^* V_\rho \cong V_{\rho \circ f_*} \ 
\end{equation}
between the pull-back of the bundle $V_\rho$ via $f$ and the bundle $V_{\rho \circ f_*} \to W$, where $f_* :
\pi_1(W) \to \pi_1(X)$ is the map induced by $f$ on the fundamental groups. 
\end{lemma}
{\em Proof:}
We have a commutative diagram
\begin{equation*}
\begin{split}
	\begin{diagram}
	 \node{\widetilde{W}} \arrow{e,t}{\widetilde{f}} \arrow{s} \node[1]{\widetilde{X}} \arrow{s} \\
	 \node{W} \arrow{e,t}{f} \node[1]{X,} 
	\end{diagram}
\end{split}
\end{equation*}
where the vertical maps are the universal coverings, and where $\widetilde{f}$ is the unique map turning the
diagram commutative (we work in the category of pointed topological spaces here). It is elementary to check
that the map $\widetilde{f}$ is equivariant with respect to the action of $\pi_1(W)$, where this group acts
on $\widetilde{X}$ via $f_* : \pi_1(W) \to \pi_1(X)$ and the deck transformation group of $\widetilde{X}$.
The claimed isomorphism follows then from naturality of the associated bundle construction. \qed

\begin{prop}\label{restrictions}
Suppose the four-manifold $X$ splits along the connected 3-ma\-ni\-fold $Y$ as $X = X_1 \cup_Y X_2$ into two
four-manifolds $X_1$ and $X_2$. Then any representation $\rho : \pi(X) \to SO(3)$ induces representations
$\rho_i : \pi_1(X_i) \to SO(3)$ via $\rho \circ (j_{i})_*$ where the map $j_i : X_i \hookrightarrow X$ is
the inclusion. For these representations we have
\begin{equation}\label{restriction}
	\left. V_\rho \right|_{X_i} = V_{\rho_i} \ .
\end{equation}
Conversely, given representations $\rho_i : \pi_1(X_i) \to SO(3)$ such that $\rho_1 \circ (k_1)_* = \rho_2
\circ (k_2)_* : \pi_1(Y) \to SO(3)$, where $k_i: Y \hookrightarrow X_i$ denotes the inclusion, there is a
representation $\rho: \pi_1(X) \to SO(3)$ inducing $\rho_1$ and $\rho_2$ via the respective restrictions. 
\end{prop}
{\em Proof:} This follows from the Theorem of Seifert and van Kampen and the lemma above or, equivalently, by gluing connections. \qed

\subsection{Vanishing results for Casson-type moduli spaces}
\begin{prop}\label{hurewicz}
Let $X$ be a four-manifold with $b_2^+(X) = 0$, and let $w \in H^2(X;\Z/2)$ be $\sum e_i \ (mod\ 2)$.
Suppose there is a representation $\rho : \pi_1(X) \to SO(3)$ with fixed second Stiefel-Whitney class $w$.
Then none of the Poincar\'e dual of the basis elements $e_i$ is in the image of the Hurewicz-homomorphism
$h: \pi_2(X) \to H_2(X;\Z)$.
\end{prop}
{\em Proof \cite[p. 1718]{T}:}  Suppose we have a map $f: S^2 \to X$ such that $PD(e_i) = f_* [S^2]$, where $[S^2]
\in H_2(S^2;\Z)$ denotes the fundamental cycle of $S^2$, and $PD(e_i)$ denotes the Poincar\'e dual of $e_i$.
Then we have
\begin{equation}\label{-1}
\langle w , f_*[S^2] \rangle  \equiv \langle \sum e_j , PD(e_i) \rangle = e_i^2 =  -1 \ (mod \ 2) .
\end{equation}
On the other hand, by naturality of the cohomology-homology pairing, we get
\begin{equation}\label{0}
\begin{split}
\langle w, f_*[S^2] \rangle & = \langle w_2(V_\rho) , f_*[S^2] \rangle = \langle f^* w_2(V_\rho) , [S^2]
\rangle \ .
\end{split}
\end{equation}
But the above Lemma \ref{naturality} implies that $f^*w_2(V_\rho) = w_2(f^* V_\rho) = w_2(V_{\rho \circ
f_*})$. As $S^2$ has trivial fundamental group the bundle $V_{\rho \circ f_*}$ is clearly the trivial
bundle, so the left hand side of equation (\ref{0}) must be zero modulo 2, a contradiction to equation (\ref{-1}). 
\qed
\begin{remark}
By Hopf's theorem on the cokernel of the Hurewicz-homomorphism, expressed in the exact sequence
\begin{equation*}
\pi_2(X) \to H_2(X;\Z) \to H_2(\pi_1(X);\Z) \to 0 \ , 
\end{equation*}
the fundamental group has to have non-trivial second homology in order to obtain a non-empty Casson-type
moduli space.
\end{remark}

This proposition gives a topological significance of the zero-energy instantons: If the moduli space is
non-empty then the elements $PD(e_i)$ are not representable by spheres! One might wonder whether there
exists any four-manifold where the elements $PD(e_i)$ are not representable by spheres. Certainly this
cannot be a simply connected four-manifold because of the Hurewicz-isomorphism theorem. Interestingly, the
answer is affirmative. Generalising Mumford's fake projective plane \cite{M}, Prasad and Yeung have constructed manifolds with the rational cohomology of the
complex projective space $\mathbb{CP}^2$ whose universal cover is the unit ball in $\C^2$ \cite{PY}. Such a manifold
$W$ is therefore an Eilenberg-MacLane space $K(\pi_1(W),1)$. 

Now let $Z$ be the four-manifold that we obtain from the connected sum of four $\overline{W}$, where
we do again a connected sum of the last summand with the first. The so obtained ``4-$\overline{W}$-ring'' is
diffeomorphic to 
\[
Z := \overline{W} \# \overline{W} \# \overline{W} \# \overline{W} \, \# \, S^1 \times S^3 \ =:4\overline{W}
\, \# \, S^1\times S^3 .
\]
This manifold has negative definite intersection form and has Betti-numbers $b_1(Z) = 1$ and $b_2(Z) = 4$.
In addition, no element of $H_2(Z,\Z)$ is representable by a 2-sphere, so we get no
obstruction to non-emptiness from Proposition \ref{hurewicz}.
Thus the four-manifold $Z$ is a prototype of a four-manifold on which to consider the moduli space of
$PU(2)$ instantons associated to the bundle $E \to Z$ with $c_1(E) = \sum e_i$ and $c_2(E) = - \frac{1}{4}
b_2(X)$ (and therefore of representations of $\pi_1(X) \to SO(3)$ with fixed Stiefel-Whitney class $w=\sum
e_i \ (mod\ 2)$). However, as we will see, there are no such instantons.

\begin{theorem}
\label{doldwhitneyapplication}
Let $X$ be a smooth closed negative definite four-manifold. If there is a representation $\rho :
\pi_1(X) \to SO(3)$ with Stiefel-Whitney class $w := \sum e_i \ (mod \ 2)$, then the second Betti-number
$b_2(X)$ must be divisible by four. 
\end{theorem}
{\em Proof:} 
The bundle $V_\rho$ has $w_2(V_\rho) = w$ and vanishing first Pontryagin-class $p_1(V_\rho) = 0$ because
this bundle admits a flat connection. Now the Dold-Whitney theorem \cite{DW} states that the second Stiefel-Whitney
class $w_2$ and the first Pontryagin class $p_1$ of any oriented real rank-3 bundle satisfy the equation
\begin{equation*}
	\text{P-Sq}(w_2) = p_1  \ \ (mod \ 4) \ .
\end{equation*}
Here $\text{P-Sq}: H^2(X;\Z/2) \to H^4(X;\Z/4)$ denotes the Pontryagin square, a lift of the cup-product
squaring $H^2(X;\Z/2) \to H^4(X;\Z/2)$ to the coefficient group $\Z/4$. If the class $v \in H^2(X;\Z/2)$ is
the mod-2 reduction of an integral class $c \in H^2(X;\Z)$ then the Pontryagin square is simply the mod-4
reduction of the square of $c$, i.e. 
\[
\text{P-Sq}(v) = c^2 \ \ (mod \ 4) \ .
\]
In our case the Dold-Whitney theorem thus implies that 
\[
0 = \text{P-Sq}(w) = \sum e_i^2 = - b_2(X) \ \ (mod \ 4) \ .
\]
\qed
Hence we obtain the following
\begin{theorem}
Let $X$ be a four-manifold with negative definite intersection form and suppose it admits a connected sum
decomposition $X_1 \# X_2$. Suppose $\rho: \pi_1(X) \to SO(3)$ is a representation of the fundamental
group of $X$ with fixed Stiefel-Whitney class $w = \sum
e_i \ (\text{mod } 2)$. Then both $b_2(X_1)$ and $b_2(X_2)$ must be divisible by four.
\end{theorem}

{\em Proof:}
Note first that the intersection form of both $X_1$ and $X_2$ must be diagonal. This follows from Eichler's
theorem on unique decomposition of symmetric definite forms over $\Z$, see \cite{HM}. Therefore the basis
vectors $\{e_i\}$ of $H^2(X;\Z)$ are simply given by the union of basis vectors $\{f_i\}$
of $H^2(X_1;\Z)$, diagonalising the intersection form of $X_1$, and basis vectors
$\{g_i\}$ of $H^2(X_2;\Z)$, diagonalising the intersection form of $X_2$. 

Note that $\pi_1(X_i \setminus B^4) \cong \pi_1(X_i)$. The above Proposition \ref{restrictions} now applies
yielding representations $\rho_i: \pi_1(X_i) \to SO(3)$. Its second Stiefel-Whitney class computes, using
the above equation (\ref{restriction}),
\begin{equation*}
\begin{split}
 	w_2(V_{\rho_1}) = \left. w_2(V_\rho) \right|_{X_1 \setminus B^4} = \sum f_i \ (mod \ 2) \ ,
\end{split}
\end{equation*}
and likewise for $w_2(V_{\rho_2})$. The above theorem therefore concludes the proof.
\qed

\begin{corollary}
This implies that the above considered manifold $Z = 4 \overline{W} \# S^1 \times S^3$ does not admit a
representation $\rho: \pi_1(X) \to SO(3)$ with Stiefel-Whitney class being the mod-2 reduction of the sum of
basis elements diagonalising the intersection form.
\end{corollary}
\begin{remark}
As a `converse' to the above vanishing theorem, suppose we are given a connected sum $X = X_1 \# X_2$ and representations $\rho_i : \pi_1(X_i) \to SO(3)$
with the desired Stiefel-Whitney classes on $X_i, \, i=1,2$. According to Proposition \ref{restrictions}, we obtain the
representation $\rho = \rho_1 * \rho_2 : \pi_1(X) \to SO(3)$ which has the desired Stiefel-Whitney class.  This is in contrast to well-known vanishing theorems for connected sums of manifolds with $b_2^+(X_i) > 0, i=1,2$ as in \cite[Theorem 9.3.4 and, in particular, Proposition 9.3.7]{DK}.
\end{remark}

\section{Constructing $4$-manifolds with non-empty Casson-type moduli space}

There is much interest in the relationship between the fundamental group of a $4$-manifold and its
intersection form.  The Casson-type invariant considered in this paper gives rise to the natural question of
whether there exists \emph{any} $4$-manifold $X$ with non-empty Casson-type moduli space.  In this section we
describe a construction that we hope will provide the first examples of such manifolds, by indicating how
to construct non-empty representation spaces $\mathscr{R}_{w}(\pi_1(X);SO(3))$. 

\subsection{Immersed $2$-links and negative-definite $4$-manifolds}

Let $\tilde{L} = \coprod_m S^2 \rightarrow S^4$ be a smooth immersion of $m$ $2$-spheres such that any points of self-intersection of
$\tilde{L}$ occur with negative sign and between two branches of the same component of $\tilde{L}$.  Suppose
there are $n$ self-intersections. Blowing up $n$ times and taking the proper transform we obtain an
$m$-component embedded link $L = \coprod_m S^2 \hookrightarrow \#^n \overline{\mathbb{CP}}^2$.

Each component of $L$ intersects each exceptional sphere of $\#^n\overline{\mathbb{CP}}^2$ either at no
points or at one point positively and at one point negatively (this is because each intersection point of
$\tilde{L}$ occurred within a single component and with negative sign).  Hence each component of $L$ is
trivial homologically and so the embedding of $L$ extends to a $D^2$-neighbourhood.

We do surgery on $L$ by removing $L \times D^2$ and gluing in $\coprod_m D^3 \times S^1$.  Call the
resulting $4$-manifold $X$.  The construction of $X$ was suggested by Kim Fr\o yshov.  It turns out to be very suited to our purposes; we have

\begin{lemma}

\begin{enumerate}

\item $H_1(X;\mathbb{Z}) = \oplus_m \mathbb{Z}$. 
\item $H_2(X;\mathbb{Z}) = \mathop{\oplus}_n \mathbb{Z}$.
\item There is a basis for $H_2(X;\mathbb{Z})$ with each element represented by an embedded torus $T^2 \hookrightarrow
X$.
\item The intersection form of $X$ is diagonal and negative definite.

\end{enumerate}
\end{lemma}

{\em Proof:}

Let $Y = \#^n \overline{\mathbb{CP}}^2 \setminus (L \times D^2)$ be the complement of the link $L$. 
Then

\begin{eqnarray*}
\#^n \overline{\mathbb{CP}}^2 &=& Y \cup m \, 2\rm{-handles} \cup m \, 4 \rm{-handles} \rm{,} \\
X &=& Y \cup m \, 3\rm{-handles} \cup m \, 4 \rm{-handles} \rm{.}
\end{eqnarray*}

\noindent Hence 

\begin{itemize}
\item $\chi(\#^n \overline{\mathbb{CP}}^2) - \chi(X) = 2m$
\item $H_1(X;\mathbb{Z}) = H_1(Y; \mathbb{Z})$
\item $H_1(Y; \mathbb{Z}) \subseteq \oplus_m \mathbb{Z}$ since $\#^n \overline{\mathbb{CP}}^2$ is simply connected.
\end{itemize}

\begin{figure}
\centerline{
{
\psfrag{exceptional_sphere}{Exceptional sphere $E$}
\psfrag{path_on_L}{Path on $L$}
\psfrag{bdy_of_ngh_of_L}{Boundary of neighbourhood of path on $L$}
\includegraphics[height=3in,width=4in]{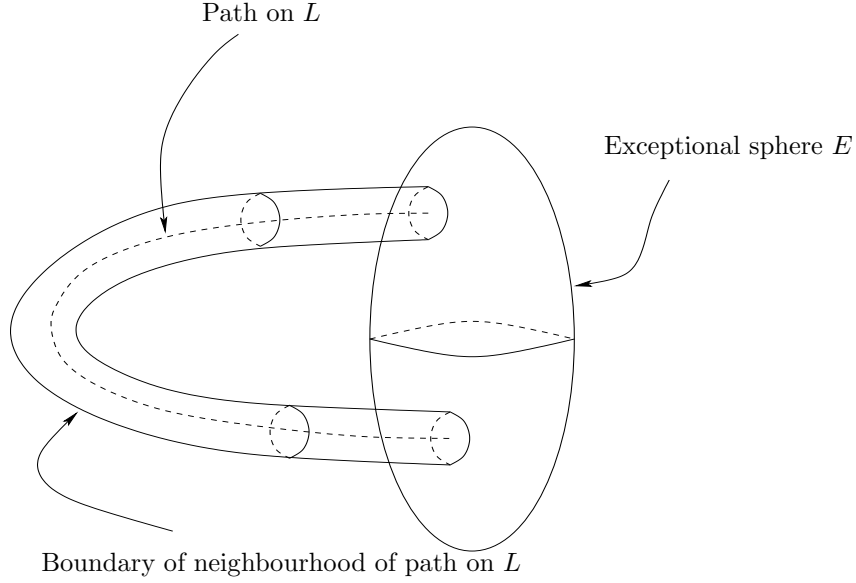}
}}
\caption{A torus representing a basis element of $H_2(X;\mathbb{Z})$.}
\label{torusbasis}
\end{figure}

\noindent So we shall be done if we can find $n$ embedded tori in $X$ which are pairwise disjoint and which
each have self-intersection $-1$.  Figure \ref{torusbasis} shows how to find these tori.  Working inside
$\#^n \overline{\mathbb{CP}}^2$, each exceptional sphere $E$ intersects $L$ transversely in two points. 
Connect these two points by a path on $L$.  The $D^2$-neighbourhood of $L$ pulls back to a trivial
$D^2$-bundle over the path.  The fibres over the two endpoints can be identified with neighbourhoods of these two points in $E$.
Removing these neighbourhoods from $E$ we get a sphere with two discs removed and we take the union of this
with the $S^1$ boundaries of all the fibres of the $D^2$-bundle over the path.

This gives a torus which has self-intersection $-1$, and we can certainly choose paths on $L$ for each
exceptional sphere which are disjoint.

\qed

We have shown how to associate to a given immersed $2$-link $\tilde{L}
= \coprod_m S^2 \rightarrow S^4$ with only negative self-intersections and disjoint components, a smooth $4$-manifold $X_{\tilde{L}}$ which is diagonal and negative
definite, with basis elements of $H_2(X_{\tilde{L}}; \mathbb{Z})$ represented by embedded tori.

\subsection{$SO(3)$ representations of $\pi_1$ and presentations of $2$-links}

Using the same notation as in the previous subsection, we give a method to describe links
$\tilde{L}$ that come with representations $\pi_1(X_{\tilde{L}}) \rightarrow SO(3)$ with the correct
Stiefel-Whitney class $w_2 = \sum e_i \, (\text{ mod } 2)$. This method may not at first appear entirely
general, but we show that if there is such a link $\tilde{L}$ then it must admit a description of this form.

\begin{figure}
\centerline{
{
\psfrag{x1}{$x_1$}
\psfrag{x4k}{$x_{4k}$}
\psfrag{y1}{$y_1$}
\psfrag{yl}{$y_l$}
\psfrag{ldots}{$\ldots$}
\includegraphics[height=3in,width=4in]{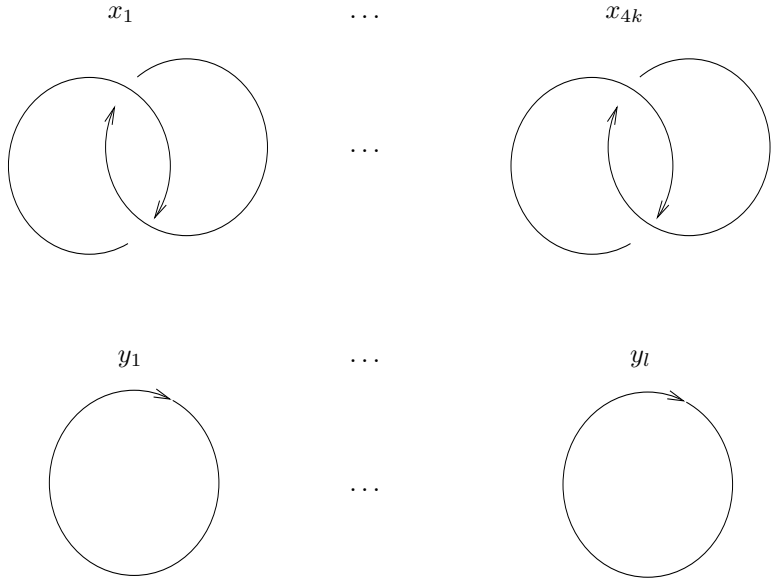}
}}
\caption{The $0$-handles of a 2-knot with $4k$ negative self-intersections.}
\label{0-handles}
\end{figure}

%%% CUT AND PASTE TO HERE.

We start by giving a lemma, which follows from basic relative Morse theory:

\begin{lemma}
\label{pres_lem}
Any closed immersed surface in $S^4$ admits a movie description in which the movie moves occur in the following
order:

\begin{enumerate}
\item $0$-handles (circle creation).
\item Simple crossing changes (see Figure \ref{crossing}).
\item Ribbon-type Reidemeister moves of type II (see Figure \ref{ribbon}).
\item Ribbon-type $1$-handle addition (see Figure \ref{ribbon}).
\end{enumerate}

\begin{figure}
\centerline{
{
\psfrag{sprout}{$0\rm{-handles} \,\, \rm{sprout} \,\, \rm{ribbons}$}
\psfrag{cross}{$\begin{array}{c}\rm{Ribbons} \,\, \rm{are} \,\, \rm{allowed} \\ \rm{to} \,\, \rm{overcross}
\,\, \rm{or} \,\, \rm{undercross} \\ \rm{each} \,\, \rm{other} \, \, \rm{and} \\ 0\rm{-handles}
\end{array}$}
\psfrag{1handle}{$\rm{Ribbon-type} \,\, 1\rm{-handle} \,\, \rm{addition}$}
\psfrag{yl}{$y_l$}
\psfrag{ldots}{$\ldots$}
\includegraphics[height=2in,width=2.5in]{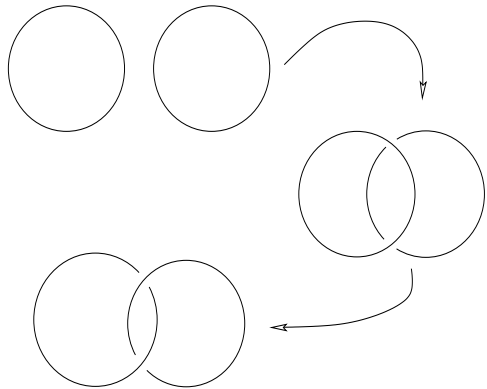}
}}
\caption{By \emph{simple crossing change} we mean doing a Reidemeister $2$ move between two $0$-crossing
diagrams of the unknot and then performing a crossing change at one of the crossings we have introduced.}
\label{crossing}
\end{figure}

After the ribbon-type $1$-handle additions there remains a diagram of an unlink and the only handle
attachments left to do are $2$-handle attachments (circle annihilation). \qed
\end{lemma}

\subsubsection{Representations of the fundamental group and ribbon presentations of $2$-links}

We now explain how to describe a representation $\pi_1 ( \#^n \overline{\mathbb{CP}}^2 \setminus L) \rightarrow SO(3)$ from a decorated presentation of the immersed link $\tilde{L}$.

For notation, let $\tilde{h} : S^4 \rightarrow \mathbb{R}$ be a height function corresponding to Lemma
\ref{pres_lem} with exactly $2$ critical points that restricts to a Morse function on $\tilde{L}$ such that
all the $i$-handles of $\tilde{L}$ occur in $\tilde{h}^{-1} ( -i )$ and the self-intersections of $\tilde{L}$ occur
in $\tilde{h}^{-1} ( -1/2 )$.  Then there is a Morse function on the blow-up $h : \#^n
\overline{\mathbb{CP}}^2 \rightarrow \mathbb{R}$ with one maximum and one minimum, and $n$ index $2$
critical points.  These index 2 critical points all occur at $h^{-1} (- 1/2)$, and h restricts to a Morse
function on the proper transform $L$ with the $i$-handles of $L$ occurring in $h^{-1} ( -i )$.  We can use the same movie of $\tilde{L}$ to describe the embedding of $L$.

Recall that $\pi_1(X_{\tilde{L}}) = \pi_1 (\#^n \overline{\mathbb{CP}}^2 \setminus L)$.  We compute
$\pi_1 (\#^n \overline{\mathbb{CP}}^2 \setminus L)$ using the Van Kampen theorem.  First note that
$h^{-1}([-3/4, \infty )) \setminus L$ is the boundary connect sum of $n$ copies of the complement of $2$ fibres in the
$D^2$-bundle over $S^2$ of Euler class $-1$, and $l$ copies of $D^4 \setminus D^2$ where the $D^2$ with $\partial D^2 \subset \partial D^4$ is
trivially embedded.  Here $n$ is the number of self-intersections of $\tilde{L}$ (and hence the number of
blow-ups required on the way to constructing $X_{\tilde{L}}$) and $l$ is the number of extra $0$-handles
used in the movie presentation of $\tilde{L}$ satisfying Lemma \ref{pres_lem}.  Since by assumption
$X_{\tilde{L}}$ has a non-empty Casson-type moduli space and dim$H_2(X_{\tilde{L}} ; \mathbb{Z}) = n$, we can write $n = 4k$ by Theorem \ref{doldwhitneyapplication}.

The boundary of $h^{-1}([-3/4, \infty )) \setminus L$ is shown as the complement of the link in Figure
\ref{0-handles}, with a point at infinity which we fix as the basepoint.

\begin{figure}
\centerline{
{
\psfrag{sprout}{$0\rm{-handles} \,\, \rm{sprout} \,\, \rm{ribbons}$}
\psfrag{cross}{$\begin{array}{c}\rm{Ribbons} \,\, \rm{are} \,\, \rm{allowed} \\ \rm{to} \,\, \rm{overcross}
\,\, \rm{or} \,\, \rm{undercross} \\ \rm{each} \,\, \rm{other} \, \, \rm{and} \\ 0\rm{-handles}
\end{array}$}
\psfrag{1handle}{$\rm{Ribbon-type} \,\, 1\rm{-handle} \,\, \rm{addition}$}
\psfrag{yl}{$y_l$}
\psfrag{ldots}{$\ldots$}
\includegraphics[height=4in,width=4in]{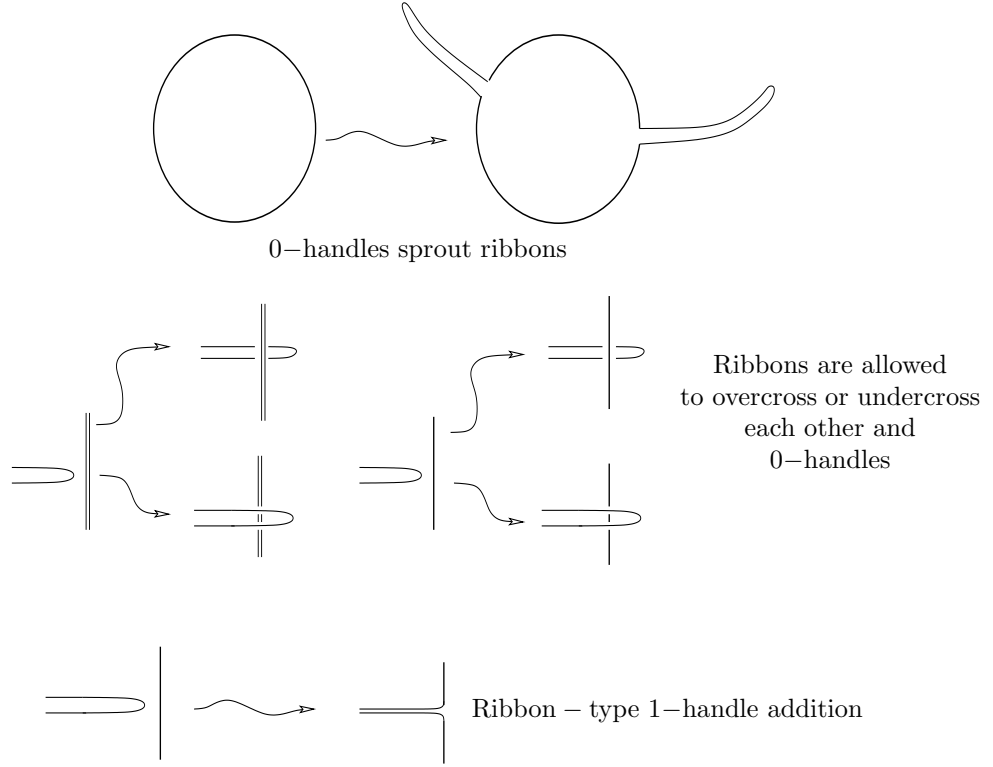}
}}
\caption{Ribbon-type moves in a movie presentation of an embedded surface in $4$-space.}
\label{ribbon}
\end{figure}

It is easy to compute that $\pi_1(h^{-1}([-3/4, \infty )) \setminus L)$ is the free
(non-abelian) group on $4k+l$ generators.  We fix representatives of a basis for this group as simple loops
coming down from infinity, linking the relevant circle by small meridians and heading back up again.  For
each of the $4k$ generators coming from the blowups we allow ourselves two representatives - one for each
circle.  Note that our representatives live in the \emph{boundary} of $h^{-1}([-3/4, \infty )) \setminus L$.

To get the space $h^{-1}([-3/2, \infty )) \setminus L$ we attach the complements of some $1$-handles to
$h^{-1}([-3/4, \infty )) \setminus L $.  What this means is that for every $1$-handle of $L$, we glue a $D^4
\setminus D^2$ to $h^{-1}([-3/4, \infty )) \setminus L$, via a homeomorphism of $(D^3\setminus (D^1 \cup
D^1)) \subseteq (S^3 \setminus S^1) = \partial(D^4 \setminus D^2)$ with a subset of $\partial (h^{-1}(-3/4)
\setminus L)$.  (All discs in this discussion are trivially embedded).

Since $\pi_1 ( D^3\setminus (D^1 \cup D^1) ) = \mathbb{Z} \times \mathbb{Z}$, $\pi_1(D^4 \setminus D^2) =
\mathbb{Z}$, and the map on $\pi_1$ induced by inclusion is onto, the Van Kampen theorem tells us that
adding the complement of a $1$-handle adds a single, possibly trivial, relation to $\pi_1$.  In other words,
we obtain a presentation of $\pi_1 ( h^{-1}([-3/2, \infty )) \setminus L )$ with $4k+l$ generators and as
many relators as there are $1$-handles.

Since we obtain $\#^n \overline{\mathbb{CP}}^2 \setminus L$ from $h^{-1}([-3/2, \infty )) \setminus L$
by gluing on the complement of some trivially embedded $D^2$'s (one for each $2$-handle of $L$) in $D^4$, it
follows that $\pi_1 ( \#^n \overline{\mathbb{CP}}^2 \setminus L ) = \pi_1 ( h^{-1}([-3/2, \infty ))
\setminus L )$.  Hence we have a presentation of $\pi_1 ( \#^n \overline{\mathbb{CP}}^2 \setminus L )$. 
Now by assumption, $X_{\tilde{L}}$ has a non-empty Casson-type moduli space, so we choose some
representation $\rho : \pi_1 ( \#^n \overline{\mathbb{CP}}^2 \setminus L ) = \pi_1(X_{\tilde{L}})
\rightarrow SO(3)$ that has the correct associated characteristic classes.  Each generator of the
presentation is associated to some circle or Hopf link in Figure \ref{0-handles}.  We decorate each circle
or Hopf link with the image of the associated generator under $\rho$.  We call these images $x_1, x_2, \ldots, x_{4k}, y_1, y_2, \ldots, y_l \in SO(3)$.

Each $1$-handle complement that we attach appears in the movie of $\tilde{L}$ as a ribbon-type $1$-handle
addition as illustrated in Figure \ref{ribbon}.  Once we have added each ribbon-type handle then by assumption we have an unlink.

%%%% %%%%This is where the CUTANDPASTE ended.

\subsubsection{Representations of the fundamental group and a singular link diagram}

We now reformulate the existence of $\rho: \pi_1 (X_{\tilde{L}}) \rightarrow SO(3)$ in terms of properties of the movie description of $\tilde{L}$ and the decoration by $x_1, \ldots, x_{4k}, y_1, \ldots, y_l \in SO(3)$.

\begin{definition}
\label{G}
A singular link diagram $G$ is given by
\begin{itemize}
\item starting with the link diagram Figure~\ref{0-handles}
\item adding the cores of each $1$-handle of $\tilde{L}$\rm{.}
\end{itemize}
\end{definition}

\noindent (For an example see Figure \ref{example}).

\begin{remark}
We could recover the full immersion $\tilde{L}$ from $G$ by adding a framing to each $1$-handle core in $G$, describing how to thicken the cores to the full $1$-handles.
\end{remark}

\begin{lemma}
\label{genus0}
These two statements are equivalent:
\begin{enumerate}
\item Each component of $\tilde{L}$ has genus $0$.
\item Suppose two circles of Figure~\ref{0-handles} are joined by three paths of $1$-handle cores $l_1, l_2, l_3$ in the singular diagram $G$.  If $l_1, l_2, l_3$ meet the first circle in three points that go clockwise (respectively anticlockwise) around the circle, then $l_1, l_2, l_3$ must meet the second circle in three points that go anticlockwise (respectively clockwise) around that circle.
\end{enumerate}
\end{lemma}

\begin{lemma}
\label{selfint}
These two statements are equivalent:
\begin{enumerate}
\item Self-intersections of $\tilde{L}$ only occur within a component and not between two components of the preimage of $\tilde{L}$.
\item The singular diagram $G$ describes an obvious singular link in $\R^3$.  Given a Hopf link in Figure~\ref{0-handles}, we require that the two circles comprising it are part of the same component in this singular link.
\end{enumerate}
\end{lemma}

The proofs of Lemmas \ref{genus0} and \ref{selfint} are left as an exercise.

\begin{lemma}
\label{relatorsgoto0}
These two statements are equivalent:
\begin{enumerate}
\item The representation

\[ \pi_1(h^{-1}([-3/4, \infty )) \setminus L) \rightarrow SO(3) \]

\noindent determined by the labelling $x_1, x_2, \ldots, x_{4k}, y_1, y_2, \ldots, y_l \in SO(3)$ factors through

\[  \pi_1 ( \#^n \overline{\mathbb{CP}}^2 \setminus L ) \rm{.}\]
\item Each circle in Figure \ref{0-handles} bounds an obvious oriented disc which has no double points when projected to the plane of the diagram.  Consider a core of a $1$-handle $A$ in the singular link diagram $G$.  Suppose $A$ connects circles
decorated by $SO(3)$ elements $g$ and $h$, and that the arc, given the orientation from $g$ to $h$,
intersects discs bounded by circles which are decorated by elements $g_1, g_2, \ldots, g_m$.  Define the
element $C(A) = (\prod_1^m g_i^{\pm 1})$, where the $\pm 1$ index is the sign of the intersection of the arc with the disc.  We require

\[ h = C(A) g C(A)^{-1} \rm{.} \]
\end{enumerate}
\end{lemma}

\begin{proof}
The condition that the given representation

\[ \pi_1(h^{-1}([-3/4, \infty )) \setminus L) \rightarrow SO(3) \]

\noindent factors through

\[  \pi_1 ( \#^n \overline{\mathbb{CP}}^2 \setminus L ) \]

\noindent is equivalent to the representation killing the relators (coming from each $1$-handle of $\tilde{L}$) in the presentation of $\pi_1 ( \#^n \overline{\mathbb{CP}}^2 \setminus L )$ discussed above.

The calculation of the relators is illustrated in
Figure \ref{relator}.

\begin{figure}
\centerline{
{
\psfrag{inf}{$\infty$}
\psfrag{g}{$g$}
\psfrag{tildeg}{$\tilde{g}$}
\psfrag{g1}{$g_1$}
\psfrag{g2}{$g_2$}
\includegraphics[height=2.5in,width=4in]{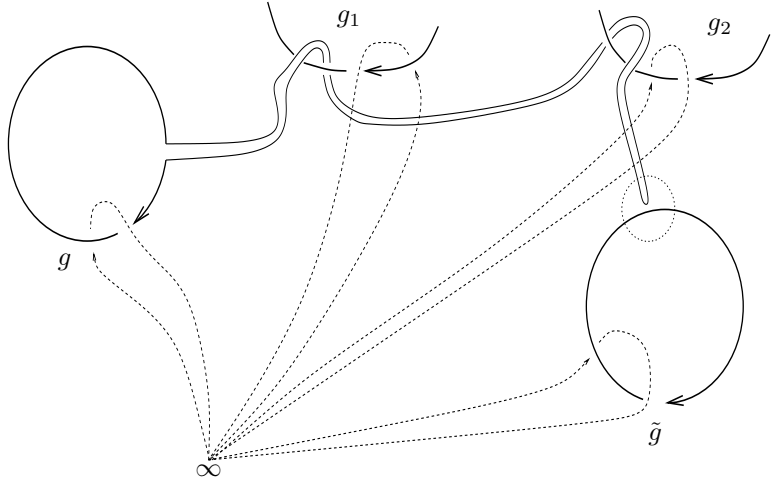}
}}
\caption{This diagram shows the situation just before the addition of a $1$-handle, which will take place
within the dotted circle.  We have indicated $4$ generators of $\pi_1 ( \#^n \overline{\mathbb{CP}}^2
\setminus L )$.  By the Van Kampen theorem, adding the $1$-handle imposes the relation that the rightmost
generator is a conjugate of the leftmost generator as in item $4$ of our checklist.  Thinking of the ribbon as a thickened arc, we note that this calculation
does not depend on whether the arc is locally knotted, but only on the order in and parity with which it
intersects the discs bounded by the $0$-handles.  Also, since a small loop encircling both strands of a
ribbon clearly bounds a disc in $h^{-1}([-3/2, \infty )) \setminus L$, it is also immaterial how the arcs
link each other.}
\label{relator}
\end{figure}
\end{proof}

\begin{lemma}
\label{swclass}
These two statements are equivalent:
\begin{enumerate}
\item The representation $\rho : \pi_1(X_{\tilde{L}}) \rightarrow SO(3)$ has the correct Stiefel-Whitney class $w_2 = \sum e_i \, (\text{ mod } 2)$.
\item \begin{itemize}
\item For each Hopf link of Figure \ref{0-handles}, choose a path of cores of $1$-handles in the singular diagram $G$ which connects the circles of the Hopf link.  (Such a path exists by Lemma \ref{selfint}).

 Say it consists of cores $A_1, A_2, \ldots, A_m$.  We order and orient these cores so that the start point of $A_1$ and the end point of $A_m$ are on different components of the Hopf link and the end point of $A_i$ is on the same circle of Figure \ref{0-handles} as the start point of $A_{i+1}$ for $1 \leq i \leq m-1$.  Write $g$ for the element decorating the Hopf link.  Then we require that

\[ \prod_1^m C(A_i) \not= 1,g \rm{.} \]
\item The elements $x_1, \ldots, x_{4k} \in SO(3)$ are each conjugate to the element $diag(1,-1,-1)$ (in other words each element $x_i$ is a rotation by $\pi$
radians).
\end{itemize}
\end{enumerate}
\end{lemma}

\begin{proof}
The condition that the representation $\rho : \pi_1(X_{\tilde{L}}) \rightarrow SO(3)$ has the correct Stiefel-Whitney class says that: 

\[ w_2(i^*\rho) = i^*(w_2(\rho)) \not= 0 \in H^2 (T^2; \mathbb{Z}/2) = \mathbb{Z}/2 \rm{,} \]

\noindent for the representative $i: T^2 \rightarrow X$ of each basis element of $H_2(X;\mathbb{Z})$.

By naturality, this means that the map $\rho \circ i : \mathbb{Z} \oplus \mathbb{Z} = \pi_1(T^2) \rightarrow
SO(3)$ has to give the non-trivial flat bundle over $T^2$.  Given a basis for $\pi_1(T^2)$ this is
equivalent to asking that $\rho \circ i$ sends each of the two basis elements to rotations by $\pi$, but around orthogonal axes.  For each basis element of $H_2(X_{\tilde{L}} ; \mathbb{Z} )$, there is an associated Hopf link in
Figure \ref{0-handles}.  Say the Hopf link is decorated by $g \in SO(3)$ and there is a path connecting the
two components of the Hopf link as in Lemma \ref{selfint}.

Consider the loop which we gave as a generator of $\pi_1 (h^{-1}([-3/4, \infty )) \setminus L)$
corresponding to the Hopf link, and a loop based at $\infty$ which goes down to the path and follows it
around until returning to the Hopf link and then returns back up to $\infty$.  This gives two basis elements
for $\pi_1$ of a $T^2$ representing the basis element of $H_2(X_{\tilde{L}} ; \mathbb{Z})$.  The former is sent to $g$ by
$\rho \circ i$ and the latter is sent to $\prod_1^m C(A_i) $.  Since necessarily $\prod_1^m C(A_i) $
commutes with $g$, the requirement that $\prod_1^m C(A_i) \not= 1,g$, ensures that  $\prod_1^m C(A_i) $ is a
rotation by $\pi$ around an axis orthogonal to that of $g$.
\end{proof}

If we can find a presentation of some $\tilde{L}$ with decoration by some $x_1, \ldots, x_{4k}, y_1, \ldots, y_l \in SO(3)$ satisfying the conditions of Lemmas \ref{genus0}, \ref{selfint}, \ref{relatorsgoto0}, and \ref{swclass}, then we have seen that we can construct a negative definite $4$-manifold $X_{\tilde{L}}$ with non-empty Casson-type moduli space.  In particular we have exhibited a particular representation

\[ \pi_1 (X_{\tilde{L}}) \rightarrow SO(3) \]

\noindent which has the required associated Stiefel-Whitney class.

Giving such a presentation of $\tilde{L}$ is equivalent to giving first the singular link diagram $G$ and then giving a framing to the cores of each $1$-handle.  Therefore we have the following:

\begin{theorem}
\label{constructionresult}
Suppose we give a singular link diagram $G$ in the sense of Definition \ref{G}, starting with Figure \ref{0-handles} and then adding arcs which begin and end at points of Figure \ref{0-handles}.  Further suppose that there is a decoration of $G$ by $x_1, \ldots, x_{4k}, y_1, \ldots, y_l \in SO(3)$ that satisfies the conditions on the singular link diagrams given as the latter statements of Lemmas \ref{genus0}, \ref{selfint}, \ref{relatorsgoto0}, and \ref{swclass}.

Then, if there exists a framing of the arcs of $G$ such that the corresponding $1$-handle additions to Figure \ref{0-handles} gives a diagram of a trivial link, there exists a $4$-manifold with non-empty Casson-type moduli space. \qed
\end{theorem}

\subsection{A partial example}

The symmetry group on $4$ elements $S_4$ can be embedded in $SO(3)$ as the rotational symmetry group of a
cube.  Under this embedding, all elements of $\rm{order}=2$ are taken to rotations by $\pi$ around some axis.

\begin{figure}
\centerline{
{
\psfrag{12}{$(12)$}
\psfrag{34}{$(34)$}
\psfrag{14}{$(14)$}
\psfrag{23}{$(23)$}
\psfrag{24}{$(24)$}
\includegraphics[height=3in,width=4in]{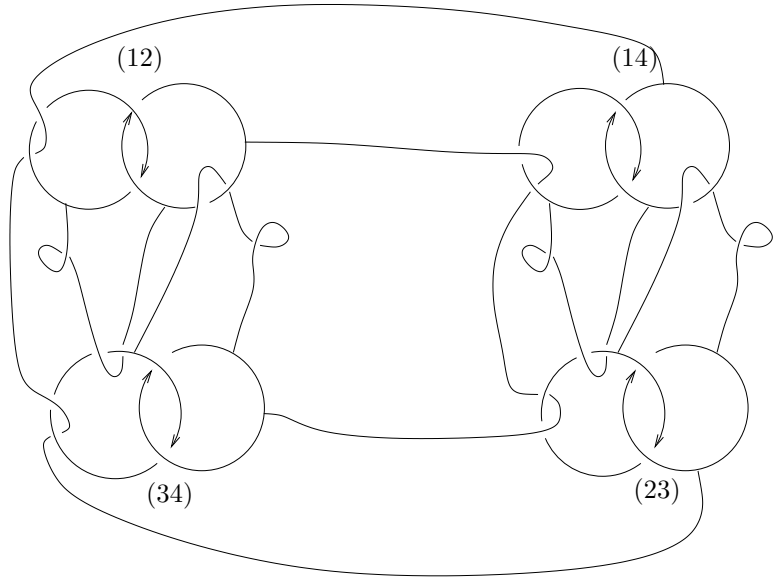}
}}
\caption{An example of what a partial diagrammatic description of a suitable immersion $\tilde{L}
\rightarrow S^4$ may look like.}
\label{example}
\end{figure}

In Figure \ref{example} we have given an example of a diagram (of labelled Hopf links, simple circles, and
arcs) satisfying all the conditions of Theorem \ref{constructionresult}.  The group element decorations of the
simple circles and the Hopf links are given in the cycle notation for $S_4 \hookrightarrow SO(3)$.  If we
can find a way to add more arcs, each satisfying the conditions of Theorem \ref{constructionresult}
(the Steifel-Whitney condition of Lemma \ref{swclass} has already been satisfied in the diagram) such that when we
replace each arc by a ribbon we get the unlink, then we will have described an
immersion $\tilde{L} \rightarrow S^4$ such that $X_{\tilde{L}}$ has non-empty Casson-type moduli space.

\begin{remark}
Note that each component of Figure \ref{example} (after replacing each arc by a blackboard-framed $1$-handle) is a smoothly slice knot.
\end{remark}

In fact in this case, more is true:

\begin{prop}
\label{onepoint}
If Figure \ref{example} is an intermediary diagram of a movie presentation of an immersed $\tilde{L}$ satisfying Theorem \ref{constructionresult}, then $X_{\tilde{L}}$ has exactly $1$ point in the representation space $\mathscr{R}_w(\pi_1(X);SO(3))$.
\end{prop}

\begin{remark}
Recent discoveries \cite{Z} have indicated that the invariant defined as the signed count of the Casson moduli space may always be $0$.  Results such as Proposition \ref{onepoint} are still valuable as they may be useful in showing that links are not slice (for more in this direction see \cite{L}).
\end{remark}

{\em Proof.}  A representation $\rho: \pi_1(X_{\tilde{L}}) \rightarrow SO(3)$ is determined by the decoration of the four Hopf links by elements of $SO(3)$.  We will see that there is only one possible decoration up to conjugation.

Suppose that we have some new decoration satisfying Theorem \ref{constructionresult}.  Call the decorating elements of $SO(3)$ $TL, TR, BL, BR$ where the initials stand for \emph{T}op, \emph{B}ottom, \emph{L}eft, \emph{R}ight.  Each of the four decorations is a rotation by $\pi$ around some axis, so each element is equivalent to a choice of axis, and we use the same labels for these axes.  By Lemma \ref{swclass}, we must have $TL$ perpendicular to $BL$ and $TR$ perpendicular to $BR$.

There is an arc connecting $TL$ to $BL$.  By condition Lemma \ref{relatorsgoto0} we can interpret this as meaning that the unique axis perpendicular to both $TR$ and $BR$ lies in the same plane as $TL$ and $BL$ and is at an angle of $\pi / 4$ to both of them.  Similarly, there is an arc connecting $TR$ and $BR$, which implies that the axis perpendicular to $TL$ and $BL$ is in the same plane as $TR$ and $BR$ and at an angle of $\pi/4$ to both of them.

It is a simple matter to convince oneself that any two ordered pairs of ordered pairs of perpendicular axes satisfying the conditions of the previous paragraph must be equivalent via the action of an element of $SO(3)$.  Hence, up to conjugation, there is exactly $1$ representation $\rho: \pi_1(X_{\tilde{L}}) \rightarrow SO(3)$ of the correct characteristic class. \qed

\end{document}